# PATH DECOMPOSITIONS FOR MARKOV CHAINS

By Götz Kersting and Kaya Memişoğlu

*University of Frankfurt*

We present two path decompositions of Markov chains (with general state space) by means of harmonic functions, which are dual to each other. They can be seen as a generalization of Williams' decomposition of a Brownian motion with drift. The results may be illustrated by a multitude of examples, but we confine ourselves to different types of random walks and the Pólya urn.

**1. Introduction and main results.** In his paper on path decomposition and local time for diffusions Williams has given an appealing decomposition of a Brownian motion with negative drift (say $-1$) at its global maximum $M$ [Williams (1974) and Pitman (1975)]. Let us state it as follows.

THEOREM 1 (Decomposition of a Brownian motion). *Let $\hat X = (\hat X_t)_{t\ge 0}$ be a Brownian motion starting in $0$ with drift $1$ and let $M$ be an independent random variable with exponential distribution and expectation $1/2$. Define the stopping time*

$$T := \sup\{t \ge 0 : \hat X_s < M \text{ for all } s < t\}.$$

*Moreover, let $\check X = (\check X_t)_{t\ge 0}$ be a process starting in $\check X_0 = \hat X_T (= M)$, whose conditional distribution, given $\hat X$ and $M$, is equal to that of a Brownian motion with drift $-1$, which is conditioned to stay below $M$. Then the process $\bar X = (\bar X_t)_{t\ge 0}$ given by*

$$\bar X_t := \begin{cases} \hat X_t, & t < T, \\ \check X_{t-T}, & t \ge T, \end{cases}$$

*is a Brownian motion with drift $-1$.*

This theorem has been the starting point for further investigations. Williams already generalized his result to one-dimensional diffusion processes, and Bertoin and Chaumont gave related paths decompositions for certain classes of Lévy processes [Bertoin (1991, 1992, 1993) and Chaumont (1996)]. Millar (1978, 1977), Jacobsen (1974) and Greenwood and Pitman (1980) discussed path decompositions from a broader point of view. Other path decompositions like Tanaka's construction arose from conditioned random walks [Tanaka (1989, 1990)]. However, these results do not seem to suggest a general pattern how to deduce path decompositions for other Markov processes.

In this paper we introduce a general method of path decomposing Markov processes by means of positive harmonic functions, which covers most of the mentioned results and thus offers a framework. We restrict ourselves to Markov chains, that is, to the case of a discrete time parameter. The technically more involved case of continuous time will be treated elsewhere.

Let $P(x,dy)$ be a probability kernel on some state space $(S,\mathcal{S})$. No (topological) restrictions are required for the state space. In the sequel $X = (X_n)_{n \in \mathbb{N}_0}$ signifies a Markov chain with transition kernel $P$ and time parameter $n \in \mathbb{N}_0 := \{0,1,2,3,\ldots\}$. The corresponding probability measure is denoted by $\mathbf{P}_x$, where $x$ as usual is the initial state of $(X_n)$.

Further, let
$$h: S \to \mathbb{R}, \qquad 0 \le h < \infty$$
be a nonnegative harmonic function with respect to $P$, that is,
$$h(x) = \int P(x,dy) h(y)$$
for all $x \in S$. Recall that for any nonnegative harmonic function $h$ we may define the $h$-transformed kernel $P^h$, given by
$$P^h(x,dy) := \frac{1}{h(x)} P(x,dy) h(y)$$
for all $x \in S$ with $h(x) > 0$. $P^h$ is a probability kernel on the restricted state space
$$S^h := \{x \in S : h(x) > 0\}.$$
Furthermore, $h$ gives rise to probability measures $\mathbf{P}_x^h$ with $x \in S^h$, given by
$$\mathbf{E}_x^h \phi(X_1,\ldots,X_n) = h(x)^{-1} \mathbf{E}_x \phi(X_1,\ldots,X_n) h(X_n) \tag{1}$$
for any measurable function $\phi: S^n \to \mathbb{R}$. As is well known, under the measure $\mathbf{P}_x^h$ the process $(X_n)$ is a Markov chain with transition kernel $P^h$.

We are ready now to state the first main result of this paper, which gives a pathwise construction of $(X_n)$ on a richer probability space. To emphasize this, we denote the corresponding probability measure by $\mathbb{P}$. Later in the section on random walks we will explain how this result fits to Williams' decomposition theorem.

THEOREM 2 (Decomposition of a Markov chain). *Let $h \geq 0$ be harmonic and $o \in S$ such that $h(o) > 0$. Let $\hat{X} = (\hat{X}_n)$ be a Markov chain with transition kernel $P^h$ and initial state $o$, defined on a probability space with probability measure $\mathbb{P}$, and let $Y$ be an independent random variable with values in $(h(o), \infty)$ and distribution given by*

$$\mathbb{P}(Y > y) = \frac{h(o)}{y},$$

*such that $Y^{-1}$ is uniformly distributed on $(0, h(o)^{-1})$. Define the random variable*

$$T := \sup\{n \geq 0 : h(\hat{X}_m) < h(\hat{X}_n) \leq Y \text{ for all } m < n\}.$$

*Moreover, let $\check{X} = (\check{X}_n)$ be a process starting in $\check{X}_0 = \hat{X}_T$, whose conditional distribution, given $\hat{X}$ and $Y$, is equal to that of a Markov chain with transition kernel $P$, which is conditioned to stay inside $\{x \in S : h(x) \leq h(\hat{X}_T)\}$. Then the process $\bar{X} = (\bar{X}_n)$ given by*

$$\bar{X}_n := \begin{cases} \hat{X}_n, & n < T, \\ \check{X}_{n-T}, & n \geq T, \end{cases}$$

*is a Markov chain with transition kernel $P$, that is, equal in distribution to $(X_n)$ under the measure $\mathbf{P}_o$.*

Variable $T$ is the moment, when $h(\bar{X}_n)$ attains its global maximum for the first time. There are cases when $T$ may also take the value $\infty$ with positive probability (in Section 3 we discuss this possibility in more detail), then no global maximum exists. In this case Theorem 2 holds with $\bar{X}_n = \hat{X}_n$ for all $n$—this simply means that no concatenation with a process $\check{X}$ takes place in this case.

Note that $T$, in general, is not a stopping time, neither for $\bar{X}$ nor for $\hat{X}$: It is the last moment, when $h(\hat{X}_n)$ exceeds all previous values, before $h(\hat{X}_n)$ surpasses $Y$ for the first time.

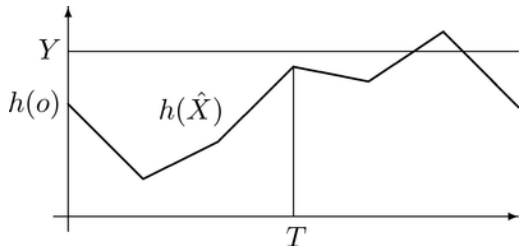

Variable $T$, thus, may contain information about the future behavior of $\hat{X}$. The notable fact is that $T$ is a *splitting time*, for $\bar{X}$ as well as for $\hat{X}$. Splitting times for Markov processes have been introduced by Jacobsen (1974) just in the context of path decompositions, they fulfil a generalized Markov property, which corresponds to the Markovian character of $\bar X$ after the moment $T$. In contrast to the classic definition of a stopping time, a splitting time $T$ allows a change of law of the process $(\bar X_{T+n})_{n\in\mathbb{N}}$, depending on the value of $T$ and $\bar X_T$, as in our case. For details compare the original work of Jacobsen (1974).

It was remarked by Doob and others that there is a kind of duality between the kernels $P$ and $P^h$ [see Chapter 12.4 in Dellacherie and Meyer (1988)]. This is also reflected in our context: The next result presents a path decomposition for a Markov chain $X = (X_n)$, whose transitions now obey the $h$-transformed kernel $P^h$. This dual decomposition takes place at the minimum of $(h(X_n))$.

THEOREM 3 (Dual decomposition). *Let $h \geq 0$ be harmonic (for the kernel $P$) and $o \in S$ such that $h(o) > 0$. Let $\check X^* = (\check X_n^*)$ be a Markov chain with transition kernel $P$ and initial state $o$ and let $U$ be an independent random variable, uniformly distributed in the interval $(0, h(o))$. Define the random variable*

$$T^* := \sup\{n \geq 0 : h(\check X_m^*) > h(\check X_n^*) \geq U \text{ for all } m < n\}.$$

*Moreover, let $\hat X^* = (\hat X_n^*)$ be a process starting in $\hat X_0^* = \check X_{T^*}^*$, whose conditional distribution, given $\check X^*$ and $U$, is equal to that of a Markov chain with transition kernel $P^h$, which is conditioned to stay inside $\{x \in S : h(x) \geq h(\check X_{T^*}^*)\}$. Then the process $\bar X^* = (\bar X_n^*)$ given by*

$$\bar X_n^* := \begin{cases} \check X_n^*, & n < T^*, \\ \hat X_{n-T}^*, & n \geq T^*, \end{cases}$$

*is a Markov chain with transition kernel $P^h$, that is, equal in distribution to $(X_n)$ under the measure $\mathbf{P}_o^h$.*

REMARK 1 [Exact sampling of $\sup_n h(X_n)$]. $\sup_n h(X_n)$ cannot be simulated directly because this would require infinitely many values $X_n$. Given that $T$ can be determined by finitely many simulation steps, our theorem offers an alternative: $\sup_n h(X_n)$ and $h(\hat X_T)$ are equal in distribution, and the theorem gives a recipe to sample the latter random variable exactly. In Section 3 we give an additional condition that ensures that $T$ can be determined algorithmicly.

REMARK 2. The conditioned Markov chains $\check X$, respectively $\hat X^*$, can also be obtained as unconditioned Markov chains, as we shall explain in the next section. In general, they are $h$-processes of $X$ only, given the value of $(\hat X_T, Y)$, respectively, $(\check X_{T^*}^*, U)$ (see Lemma 5 for details).

REMARK 3. If $g = ch$ for some constant $c > 0$, then $P^h$ and $P^g$ are equal. Correspondingly, our constructions are the same for $g$ and $h$ up to the scaling factor $c$. One would also expect that, in essence, the path decompositions remains unaltered, if $h$ is replaced by $h + c$. This is not easy to see directly for the first construction. Here duality is helpful, the transition from $h$ to $h+c$ can be well understood for the second construction. For details we refer to the proof of Theorem 3.

REMARK 4 (Doob inequality). As an immediate consequence of our decomposition, we obtain the inequality

$$\mathbb{P}\left(\sup_{1 \leq i \leq n} h(\bar{X}_i) \geq \lambda\right) \leq \mathbb{P}(Y > \lambda) = \frac{h(\bar{X}_0)}{\lambda} = \frac{\mathbb{E}h(\bar{X}_n)}{\lambda}.$$

This is Doob's martingale inequality, applied to the martingale $(h(X_n))$.

The paper is organized as follows. In Section 2 we prove both theorems. In Section 3 the question is addressed, how to detect the values of $T$ and $T^*$ algorithmicly. In Section 4 we apply the decompositions to different types of random walks.

**2. Proof of the main results.** Let us first have a closer look at the conditioned chain described in Theorem 2. It can be obtained as a Markov chain with suitable transition kernel. Let us recall this well-known construction. For $s \geq 0$ define

$$q_s(x) := \mathbf{P}_x(h(X_i) \leq s \, \forall i \in \mathbb{N}_0)$$

and the stopping time

$$\sigma_s := \inf\{i \geq 0 : h(X_i) > s\}.$$

LEMMA 4. *If $h(x) \leq s$, then $q_s(x) > 0$.*

PROOF. Assume $\sigma_s < \infty$ a.s. Then Fatou's lemma and the martingale property of $(h(X_n))$ lead to

$$\mathbf{E}_x h(X_{\sigma_s}) = \mathbf{E}_x \lim_{n \to \infty} h(X_{\sigma_s \wedge n}) \leq \liminf_{n \to \infty} \mathbf{E}_x h(X_{\sigma_s \wedge n}) = h(x) \leq s.$$

Because $h(X_{\sigma_s}) > s$ a.s., this is a contradiction. Thus, $\mathbf{P}_x(\sigma_s = \infty) = q_s(x) > 0$ and the lemma is proved. □

Furthermore, for $h(x) \leq s$ the Markovian character of the chain gives

$$q_s(x) = \int_{y \,:\, h(y) \leq s} P(x, dy) q_s(y).$$

Thus, $q_s$ is also a harmonic function for the kernel $P$ restricted to the set $S_s := \{x \in S : h(x) \leq s\}$ and again for each $s$ we can define a $q_s$-transformed kernel $Q_s$ of $P$ on the respective set $S_s$ by setting

$$Q_s(x, dy) := \frac{1}{q_s(x)} P(x, dy) q_s(y)$$
$$= \mathbf{P}_x(X_1 \in dy \mid h(X_i) \leq s \,\forall\, i \in \mathbb{N}_0).$$

It follows that

$$\mathbf{P}_x(X_1 \in dx_1, \ldots, X_n \in dx_n \mid h(X_i) \leq s \,\forall\, i \in \mathbb{N}_0)$$
$$= \frac{P(x, dx_1) \cdots P(x_{n-1}, dx_n) q_s(x_n)}{q_s(x)}$$
$$= Q_s(x, dx_1) \cdots Q_s(x_{n-1}, dx_n).$$

Altogether we obtain the following result.

LEMMA 5. *Under $\mathbf{P}_x$ the process $(X_n)$, conditioned to stay inside $S_s$, $s \geq 0$, is a Markov chain with transition kernel $Q_s(x, dy)$.*

The main step of our proof, which is based on a change of measure type argument, is contained in the next lemma.

LEMMA 6. *Let $h$ be harmonic and $0 < h(o) \leq h(x) \leq s$. Let $U$ be uniformly distributed in $(0, h(o)^{-1})$ and independent of $(X_n)$ (with respect to $\mathbf{P}_x^h$). Then*

$$\mathbf{P}_x^h\left(U \leq \frac{1}{h(X_{\sigma_s})} I_{\{\sigma_s < \infty\}}\right) = \frac{h(o)}{h(x)} \mathbf{P}_x(\sigma_s < \infty),$$

$$\mathbf{P}_x^h\left(\frac{1}{h(X_{\sigma_{h(x)}})} I_{\{\sigma_{h(x)} < \infty\}} < U \leq \frac{1}{h(x)}\right) = \frac{h(o)}{h(x)} q_{h(x)}(x).$$

PROOF. Since $(h(X_n))$ is a $\mathbf{P}_x$-martingale and $\sigma_s$ a stopping time,

$$\mathbf{P}_x(\sigma_s \leq n) = \mathbf{E}_x\left(\frac{h(X_{\sigma_s \wedge n})}{h(X_{\sigma_s \wedge n})}; \sigma_s \leq n\right) = \mathbf{E}_x\left(\frac{h(X_n)}{h(X_{\sigma_s \wedge n})}; \sigma_s \leq n\right)$$

[note that $h(X_{\sigma_s \wedge n}) > 0$ on the event $\{\sigma_s \leq n\}$]. By means of (1) we rewrite this equation as

$$\mathbf{P}_x(\sigma_s \leq n) = h(x) \mathbf{E}_x^h\left(\frac{1}{h(X_{\sigma_s \wedge n})}; \sigma_s \leq n\right)$$
$$= h(x) \mathbf{E}_x^h\left(\frac{1}{h(X_{\sigma_s})}; \sigma_s \leq n\right).$$

As $n \to \infty$,

(2) $$\mathbf{P}_x(\sigma_s < \infty) = h(x)\mathbf{E}_x^h\left(\frac{1}{h(X_{\sigma_s})}I_{\{\sigma_s < \infty\}}\right)$$

and, consequently,

(3) $$q_{h(x)}(x) = 1 - \mathbf{P}_x(\sigma_{h(x)} < \infty)$$
$$= h(x)\mathbf{E}_x^h\left(\frac{1}{h(x)} - \frac{1}{h(X_{\sigma_{h(x)}})}I_{\{\sigma_{h(x)} < \infty\}}\right).$$

Since $h(X_{\sigma_s})^{-1}, h(x)^{-1} \leq h(o)^{-1}$, by assumption,

$$\mathbf{P}_x^h\left(U \leq \frac{1}{h(X_{\sigma_s})}I_{\{\sigma_s < \infty\}} \,\Big|\, (X_n)\right) = \frac{h(o)}{h(X_{\sigma_s})}I_{\{\sigma_s < \infty\}},$$

$$\mathbf{P}_x^h\left(\frac{1}{h(X_{\sigma_{h(x)}})}I_{\{\sigma_{h(x)} < \infty\}} < U \leq \frac{1}{h(x)} \,\Big|\, (X_n)\right) = \frac{h(o)}{h(x)} - \frac{h(o)}{h(X_{\sigma_{h(x)}})}I_{\{\sigma_{h(x)} < \infty\}}.$$

Taking expectations, the claim follows from (2) and (3). $\square$

PROOF OF THEOREM 2. Define

$$\tau := \inf\{i \geq 0 : h(X_j) \leq h(X_i) \,\forall j \geq i\}$$

to be the moment, where $(h(X_n))$ attains its global maximum for the first time. First, we prove

(4) $\quad \mathbf{P}_o((X_1, \ldots, X_n) \in B, \tau = m) = \mathbb{P}((\bar{X}_1, \ldots, \bar{X}_n) \in B, T = m)$

for natural numbers $0 \leq m \leq n$. Since

$$\{(X_1, \ldots, X_n) \in B, \tau = m\}$$
$$= \{(X_1, \ldots, X_n) \in B \cap B_{m,n}, h(X_j) \leq h(X_m) \,\forall j \geq n\}$$

with

$$B_{m,n} := \{(x_1, \ldots, x_n) : h(x_1), \ldots, h(x_{m-1})$$
$$< h(x_m) \geq h(x_{m+1}), \ldots, h(x_n)\},$$

it follows that

$$\mathbf{P}_o((X_1, \ldots, X_n) \in B, \tau = m)$$
$$= \int_{B \cap B_{m,n}} P(o, dx_1) \cdots P(x_{n-1}, dx_n) q_{h(x_m)}(x_n).$$

On the other hand, since $h(\hat{X}_T) = h(\check{X}_0) \geq h(\check{X}_1), h(\check{X}_2), \ldots,$

$$\{(\bar{X}_1, \ldots, \bar{X}_n) \in B, T = m\}$$
$$= \{(\hat{X}_1, \ldots, \hat{X}_m, \check{X}_1, \ldots, \check{X}_{n-m}) \in B \cap B_{m,n}\}$$
$$\cap (\{h(\hat{X}_m) \leq Y < h(\hat{X}_{\hat{\sigma}_m}), \hat{\sigma}_m < \infty\} \cup \{h(\hat{X}_m) \leq Y, \hat{\sigma}_m = \infty\}),$$

where the stopping times $\hat{\sigma}_m$ are defined as

$$\hat{\sigma}_m := \inf\{i \geq m : h(\hat{X}_i) > h(\hat{X}_m)\}.$$

It follows that

$$\mathbb{P}((\bar{X}_1, \ldots, \bar{X}_n) \in B, T = m)$$
$$= \int_{B \cap B_{m,n}} P^h(o, dx_1) \cdots P^h(x_{m-1}, dx_m)$$
$$\times Q_{h(x_m)}(x_m, dx_{m+1}) \cdots Q_{h(x_m)}(x_{n-1}, dx_n)$$
$$\times \mathbb{P}\left(\frac{1}{h(\hat{X}_{\hat{\sigma}_m})} I_{\{\hat{\sigma}_m < \infty\}} < \frac{1}{Y} \leq \frac{1}{h(\hat{X}_m)} \Big| \hat{X}_m = x_m\right)$$
$$= \int_{B \cap B_{m,n}} P(o, dx_1) \cdots P(x_{n-1}, dx_n) \frac{h(x_m) q_{h(x_m)}(x_n)}{h(o) q_{h(x_m)}(x_m)}$$
$$\times \mathbf{P}^h_{x_m}\left(\frac{1}{h(X_{\sigma_{h(x_m)}})} I_{\{\sigma_{h(x_m)} < \infty\}} < \frac{1}{Y} \leq \frac{1}{h(x_m)}\right).$$

Now (4) follows in view of Lemma 6 because $Y^{-1}$ is uniformly distributed.

Next, we prove

(5) $\quad \mathbf{P}_o((X_1, \ldots, X_n) \in B, \tau > n) = \mathbb{P}((\bar{X}_1, \ldots, \bar{X}_n) \in B, T > n).$

From

$$\{(X_1, \ldots, X_n) \in B, \tau > n\}$$
$$= \{(X_1, \ldots, X_n) \in B\}$$
$$\cap \{h(X_i) > \max(h(X_1), \ldots, h(X_n)) \text{ for some } i > n\},$$

it follows that

$$\mathbf{P}_o((X_1, \ldots, X_n) \in B, \tau > n)$$
$$= \int_B P(o, dx_1) \cdots P(x_{n-1}, dx_n) \mathbf{P}_{x_n}(\sigma_{\max(h(x_1), \ldots, h(x_n))} < \infty).$$

On the other hand,

$$\{(\bar{X}_1, \ldots, \bar{X}_n) \in B, T > n\}$$
$$= \{(\hat{X}_1, \ldots, \hat{X}_n) \in B\} \cap \{\bar{\sigma}_n < \infty, h(\hat{X}_{\bar{\sigma}_n}) \leq Y\}$$

with
$$\overline{\sigma}_n := \inf\{i > n : h(\hat{X}_i) > \max(h(\hat{X}_1), \ldots, h(\hat{X}_n))\}.$$

It follows that
$$\mathbb{P}((\bar{X}_1, \ldots, \bar{X}_n) \in B, T > n)$$
$$= \int_B P^h(o, dx_1) \cdots P^h(x_{n-1}, dx_n)$$
$$\times \mathbb{P}\left(\frac{1}{Y} \leq \frac{1}{h(\hat{X}_{\overline{\sigma}_n})} I_{\{\overline{\sigma}_n < \infty\}} \bigg| \hat{X}_n = x_n\right)$$
$$= \int_B P(o, dx_1) \cdots P(x_{n-1}, dx_n) \frac{h(x_n)}{h(o)}$$
$$\times \mathbf{P}^h_{x_n}\left(\frac{1}{Y} \leq \frac{1}{h(X_{\sigma_{\max(h(x_1),\ldots,h(x_n))}})} I_{\{\sigma_{\max(h(x_1),\ldots,h(x_n))} < \infty\}}\right).$$

Since $Y^{-1}$ is uniformly distributed and because of Lemma 6, (5) also follows.

Combining (4) and (5), we obtain
$$\mathbb{P}((\bar{X}_1, \ldots, \bar{X}_n) \in B) = \mathbf{P}_o((X_1, \ldots, X_n) \in B)$$

and the theorem is proved. □

Theorem 3 could be proved along the same lines, with some marked differences. We prefer to deduce it from Theorem 2 via duality.

PROOF OF THEOREM 3. First assume that $h(x) > 0$ for all $x \in S$. Then, as is easily seen, $1/h(x)$ is harmonic with respect to $P^h$, moreover, the $1/h$-transformed kernel of $P^h$ is equal to $P$, $(P^h)^{1/h} = P$. Also $Y^{-1}$ is uniformly distributed in $(0, (1/h(o))^{-1}) = (0, h(o))$. In this situation Theorem 3 is an immediate consequence of Theorem 2.

Next, let us assume that $h(x)$ also takes the value 0. Then $h(x) + \varepsilon$ is a strictly positive harmonic function for $\varepsilon > 0$ with respect to $P$. The formula

$$(6) \qquad \mathbf{P}_o^{h+\varepsilon} = \frac{h(o)}{h(o) + \varepsilon} \mathbf{P}_o^h + \frac{\varepsilon}{h(o) + \varepsilon} \mathbf{P}_o,$$

which follows from (1), enables us to tranfer the validity of the path decomposition according to Theorem 3 from $h + \varepsilon$ to $h$. On the one hand $\mathbf{P}_o^{h+\varepsilon}$ converges to $\mathbf{P}_o^h$, as follows from (6). On the other hand, the processes induced by our construction (resp. its components, the uniform random variables, the splitting times, and the conditioned Markov chains) converge in distribution, as $\varepsilon \to 0$. □

**3. Detecting the value of the splitting time.** It is too much to expect that such a general decomposition result as ours will provide useful information for any Markov chain possessing nontrivial harmonic functions. There are examples with an exotic touch.

As an attempt to distinguish favorable situations from less appealing cases, let us adopt an algorithmic point of view and ask the question: When is it sufficient to run the chain $\hat X$ for finitely many steps in order to detect the value of $T$? Recall that, in general, $T$ is not a stopping time and cannot be observed immediately. Clearly $T$ is determined by $(\hat X_n)_{n \leq \tau_c}$, where the stopping time $\tau_c \geq T$ denotes the moment of crossing the level $Y$,

$$\tau_c := \inf\{n \geq 0 : h(\hat X_n) > Y\}.$$

Thus, on the event $\{\sup_n h(\hat X_n) > Y\} = \{\tau_c < \infty\}$ only finitely many steps of $\hat X$ are required.

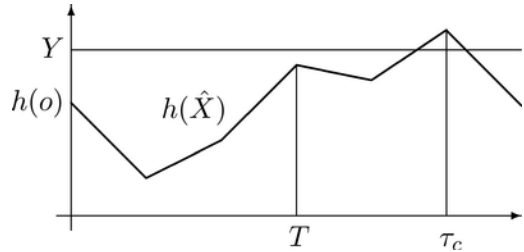

The situation may become less desirable, however, if $\{Y \geq \sup h\} = \{\tau_c = \infty\}$ occurs. Then it may happen that $T$ has a finite value, which, nevertheless, cannot be observed within finite time, neither by $\tau_c$ nor by some other mean. To put this in mathematical terms: We say that $T$ is *detected by a stopping time* $\tau$ (stopping time with respect to $\hat X$; we allow that $\tau$ depends on $Y$, too), if finite values of $T$ may be recognized from finite values of $\tau$, namely, that the value of $T$ is a.s. not bigger than $\tau$. Expressed in formulas, this means

$$\{T < \infty\} = \{\tau < \infty\} = \{T \leq \tau < \infty\}, \qquad \mathbb{P}\text{-a.s.}$$

Then $T = \sup\{n \leq \tau : n < \tau_c, h(X_n) > \max_{i<n} h(X_i)\}$, $\mathbb{P}$-a.s. Note that $\tau_c$ detects $T$, if $\{T < \infty\} = \{\tau_c < \infty\}$, $\mathbb{P}$-a.s.

An instance where no $T$-detecting stopping time exists may be considered as unfavourable and we give an example at the end of this section (see Example 1). For a majority of harmonic functions this kind of phenomenon does not occur. We introduce now a class of harmonic functions leading to path decompositions with splitting times, which always can be detected by stopping times.

First observe that the condition

(7) $$\sup_n h(X_n) = \infty, \qquad \mathbf{P}_x^h\text{-a.s.},$$

for all $x \in S^h$ implies $\sup_n h(\hat{X}_n) = \infty$ and $\tau_c < \infty$ $\mathbb{P}$-a.s. Then $T$ is detected by $\tau_c$. More generally, consider the condition

(8) $$\sup_n h(X_n) = \sup h, \qquad \mathbf{P}_x^h\text{-a.s.},$$

for all $x \in S^h$, thus, $\sup_n h(\hat{X}_n) = \sup h$ $\mathbb{P}$-a.s. Then, in case that $\{\tau_c = \infty\}$ occurs, $T$ takes the value $\infty$, if $h(\hat{X}_n) < \sup h$ for all $n$, and, otherwise, the smallest value $n$ such that $h(\hat{X}_n) = \sup h$. Thus, introducing the stopping time

$$\hat{\tau} := \inf\{n \geq 0 : h(\hat{X}_n) = \sup h\}$$

$T = \hat{\tau}$ on the event $\{\tau_c = \infty\}$. It follows that $T$ is detected by the stopping time $\tau := \tau_c \wedge \hat{\tau}$.

The most prominent class of functions fulfilling (8) is constituted by the minimal harmonic functions. Recall that a nonnegative harmonic function $h$ is said to be minimal (or extremal), if for any other harmonic function $k$ such that $0 \leq k(x) \leq h(x)$ for all $x \in S$, it follows $k = ch$ with some constant $c \in [0, 1]$.

For a countable state space the following result is a consequence of Doob's theorem from Martin boundary theory. We give a proof for arbitrary state space.

PROPOSITION 7. *If $h$ is a minimal harmonic function, then*

$$h(X_n) \to \sup h, \qquad \mathbf{P}_x^h\text{-a.s.},$$

*for all $x \in S^h$.*

PROOF. Since

$$\int_{S^h} P^h(x, dy) h(y)^{-1} = \int_{S^h} h(x)^{-1} P(x, dy) \leq h(x)^{-1} \int_S P(x, dy) = h(x)^{-1},$$

$1/h(X_n)$ is a $\mathbf{P}_x^h$-supermartingale. Thus, from the martingale convergence theorem

$$h(X_n) \to L, \qquad \mathbf{P}_x^h\text{-a.s.},$$

where $L$ is a random variable with values in $(0, \sup h]$. For $\gamma \geq 0$ define the function

$$k_\gamma(x) := \mathbf{P}_x^h(L \leq \gamma)$$

for $x \in S^h$. This function is harmonic with respect to the kernel $P^h$ because

$$k_\gamma(x) = \mathbf{P}_x^h\left(\lim_n h(X_n) \leq \gamma\right)$$

$$= \int_{S^h} P^h(x, dy) \mathbf{P}_y^h\left(\lim_n h(X_n) \leq \gamma\right) = \int_{S^h} P^h(x, dy) k_\gamma(y).$$

It follows that

$$h(x) k_\gamma(x) = \int_{S^h} P(x, dy) h(y) k_\gamma(y) = \int_S P(x, dy) h(y) k_\gamma(y),$$

where $k_\gamma(x)$ may be chosen arbitrary for $x \notin S^h$. This last equation says that $h \cdot k_\gamma$ is harmonic with respect to $P$, also

$$h \cdot k_\gamma \leq h$$

by definition of $k_\gamma$. Since by assumption $h$ is minimal, we conclude that

$$h \cdot k_\gamma = \alpha h$$

and, therefore,

$$\mathbf{P}_x^h(L \leq \gamma) = \alpha$$

for all $x \in S^h$ with some constant $\alpha = \alpha_\gamma \in [0, 1]$.

Next we show that the only possible values for $\alpha$ are 0 or 1. Consider the estimate

$$\mathbf{P}_x^h\left(\{L \leq \gamma\} \cap \bigcap_{m \geq n} \{h(X_m) \leq \gamma\}\right)$$

$$\leq \mathbf{P}_x^h(L \leq \gamma, h(X_n) \leq \gamma) \leq \mathbf{P}_x^h(L \leq \gamma).$$

Since $h(X_n) \to L$, $\mathbf{P}_x^h$-a.s., the left-hand side converges to $\mathbf{P}_x^h(L \leq \gamma)$, as $n \to \infty$, provided $\gamma$ is a point of continuity of $L$. Consequently,

$$\mathbf{P}_x^h(h(X_n) \leq \gamma, L \leq \gamma) \to \alpha.$$

On the other hand, by the Markov property,

$$\mathbf{P}_x^h(h(X_n) \leq \gamma, L \leq \gamma) = \int_{y: h(y) \leq \gamma} \mathbf{P}_x^h(h(X_n) \in dy) \mathbf{P}_y^h(L \leq \gamma)$$

$$= \int_{y: h(y) \leq \gamma} \mathbf{P}_x^h(h(X_n) \in dy) \alpha$$

$$= \alpha \mathbf{P}_x^h(h(X_n) \leq \gamma).$$

If $\gamma$ is a point of continuity of the distribution of $L$, it follows that

$$\mathbf{P}_x^h(h(X_n) \leq \gamma, L \leq \gamma) \to \alpha \mathbf{P}_x^h(L \leq \gamma) = \alpha^2.$$

Therefore, $\alpha = \alpha^2$, respectively, $\alpha = 0$ or $1$ for all $\gamma$ (up to countably many), which means that $L$ has a degenerate distribution,

$$L = \beta, \qquad \mathbf{P}_x^h\text{-a.s.}$$

for some constant $\beta \in (0, \sup h]$. Also $\beta$ does not depend on the starting state $x$ of $(X_n)$. Since $1/h(X_n)$ is a nonnegative $\mathbf{P}_x^h$-supermartingale, we conclude by Fatou's lemma,

$$h(x)^{-1} \geq \lim_n \mathbf{E}_x^h h(X_n)^{-1} \geq \mathbf{E}_x^h L^{-1} = \beta^{-1}$$

or $h(x) \leq \beta \leq \sup h$ for all $x \in S^h$. Therefore, $L = \sup h$, $\mathbf{P}_x^h$-a.s., which is our claim. $\square$

An interesting example with different harmonic functions is provided by Pólya's urn scheme. We show now that some of these functions lead to path decompositions with detectable splitting times, while others do not.

EXAMPLE 1. In *Pólya's urn scheme* each drawn ball is put back into the urn together with an extra ball of the same color (two different colors). We consider the Markov chain $X_n = (R_n, T_n)$, $n \geq 0$, with states $x = (r, t) \in \mathbb{N}^2$, $0 < r < t$, where $r$ denotes the number of red balls and $t$ the total number of balls in the urn. The transition probabilities are

$$P_{xy} = \begin{cases} \dfrac{r}{t}, & \text{for } x = (r,t),\ y = (r+1, t+1), \\ 1 - \dfrac{r}{t}, & \text{for } x = (r,t),\ y = (r, t+1). \end{cases}$$

Then for given $0 < p < 1$, $q = 1 - p$,

$$h_p(x) := (t-1) \binom{t-2}{r-1} p^{r-1} q^{t-r-1}$$

is a harmonic function. Here $\sup h_p = \infty$. A quick calculation shows $P_{xy}^{h_p} = p$, respectively, $q$, such that $\hat{X}$ is a random walk. It is not difficult to conclude that (7) is fulfilled. This also follows from the last proposition because it is well known that these harmonic functions are minimal [Blackwell and Kendall (1964)]. Consequently, the associated path decompositions possess detectable splitting times. (The corresponding dual decomposition establishes an exotic path decomposition of a random walk by means of a Pólya urn.)

Also,

$$h(x) := \frac{r}{t}$$

is a harmonic function. In this case $P_{xy}^h = \frac{r+1}{t+1}$, respectively, $1 - \frac{r+1}{t+1}$, consequently, under the measure $\mathbf{P}_x^h$ the shifted random variables $X_n + (1,1)$, $n \geq 0$, can be viewed as a Pólya urn (with one extra red ball added to the urn at the beginning). We like to show that the corresponding path decomposition gives rise to an undetectable splitting time $T$. To this end we make use of the familiar facts that $h(X_n) = R_n/T_n$ is a.s. convergent and that, given the limit $L$, the random variables $Z_n := R_n - R_{n-1}$, $n \geq 1$, are independent Bernoulli variables with success probability $L$. These statements are valid with respect to $\mathbf{P}_x$ as well as $\mathbf{P}_x^h$.

First, we show that $T < \infty$, $\mathbb{P}$-a.s. Write $\hat{X}_n = (\hat{R}_n, \hat{T}_n)$, $\hat{Z}_n = \hat{R}_n - \hat{R}_{n-1}$ and $\hat{L} = \lim_n h(\hat{X}_n)$ $\mathbb{P}$-a.s. Then from $\hat{R}_n = \hat{R}_0 + \sum_{i=1}^n \hat{Z}_i$ and $\hat{T}_n = n + \hat{T}_0$, it follows

$$\left\{\sup_n h(\hat{X}_n) \leq \hat{L}\right\} = \left\{\sup_n \sum_{i=1}^n (\hat{Z}_i - \hat{L}) \leq \hat{L}\hat{T}_0 - \hat{R}_0\right\},$$

thus, by the properties of sums of i.i.d. Bernoulli variables $\{\sup_n h(\hat{X}_n) \leq \hat{L}\}$ has zero probability. In other words, the moment

$$N := \inf\left\{n \geq 0 : h(\hat{X}_n) = \sup_m h(\hat{X}_m)\right\},$$

when for the first time $h(\hat{X}_n)$ obtains its global maximum, is finite $\mathbb{P}$-a.s., which implies $T < \infty$, $\mathbb{P}$-a.s. Furthermore, $\sup_n h(\hat{X}_n) \leq 1$, such that $\{\tau_c = \infty\} = \{\sup_n h(\hat{X}_n) \leq Y\}$ has strictly positive probability. Therefore, $\tau_c$ does not detect $T$.

In fact, there is no stopping time detecting $T$ at all. For any $\mathbb{P}$-a.s. finite stopping time $\tau$ we have $\{N \leq \tau\} = \{\sup_{j \geq 0} h(\hat{X}_{\tau+j}) \leq M\}$, $\mathbb{P}$-a.s., where $M := \max_{n \leq \tau} h(\hat{X}_n)$. Thus, from the strong Markov property

$$\mathbb{P}(N \leq \tau) = \mathbb{E}\psi(\hat{X}_\tau, M) \qquad \text{with } \psi(x, m) := \mathbf{P}_x^h\left(\sup_{j \geq 0} h(X_j) \leq m\right).$$

Since $\psi(x, m) < 1$ for $m < 1$ and since $M < 1$, $\mathbb{P}$-a.s., it follows $\mathbb{P}(N \leq \tau) < 1$. Similarly, $\mathbb{P}(N \leq \tau \mid Y) < 1$, $\mathbb{P}$-a.s. (recall that we allow $\tau$ to depend on $Y$). Because $T = N$ on the event $\{Y \geq 1\}$, we obtain

$$\mathbb{P}(Y \geq 1, T \leq \tau) = \mathbb{P}(Y \geq 1, N \leq \tau)$$
$$= \mathbb{E}(\mathbb{P}(N \leq \tau | Y); Y \geq 1) < \mathbb{P}(Y \geq 1)$$

and, consequently, $\mathbb{P}(T \leq \tau) < 1$. Since $T < \infty$, $\mathbb{P}$-a.s., $\tau$ cannot detect $T$.

REMARK 5. The class of harmonic functions given by (8) can be described in different ways. Given $x \in S^h$, we have the following equivalent conditions:

$$\sup_n h(X_n) = \sup h, \qquad \mathbf{P}_x^h\text{-a.s.}$$

$$\iff \quad \lim_n h(X_n) = \sup h, \qquad \mathbf{P}_x^h\text{-a.s.}$$

$$\iff \quad \lim_n h(X_n) \in \{0, \sup h\}, \qquad \mathbf{P}_x\text{-a.s.}$$

For the proof we refer to the following facts, which follow from the martingale property of $(h(X_n))$ with respect to $\mathbf{P}_x$ and the supermartingal property of $(1/h(X_n))$ with respect to $\mathbf{P}_x^h$: (i) $h(X_n)$ is a.s. convergent to a random variable $L$, with respect to both measures $\mathbf{P}_x$ and $\mathbf{P}_x^h$, where $L < \infty$, $\mathbf{P}_x$-a.s. and $L > 0$, $\mathbf{P}_x^h$-a.s., (ii) if $h(X_n) = \sup h$ for some $n$, then $h(X_m) = \sup h$, $\mathbf{P}_x^h$-a.s. for all $m \geq n$, (iii) the measures $\mathbf{P}_x$ and $\mathbf{P}_x^h$ are related by the formulas

$$d\mathbf{P}_x^h(\cdot \cap \{L < \infty\}) = \frac{L}{h(x)} \, d\mathbf{P}_x(\cdot), \qquad d\mathbf{P}_x(\cdot \cap \{L > 0\}) = \frac{h(x)}{L} \, d\mathbf{P}_x^h(\cdot)$$

[see Theorem 1, Section VII.6, in Shiryaev (1995)]. In much the same manner it follows that under the assumption $h(x) < \sup h$,

$$\inf_n h(X_n) = \inf h, \qquad \mathbf{P}_x\text{-a.s.}$$

$$\iff \quad \lim_n h(X_n) = \inf h, \qquad \mathbf{P}_x\text{-a.s.}$$

$$\iff \quad \lim_n h(X_n) \in \{\inf h, \infty\}, \qquad \mathbf{P}_x^h\text{-a.s.}$$

These conditions describe a class of harmonic functions such that the splitting time $T^*$ of the dual path decomposition is detectable by some stopping time. Combining both sets of conditions, we obtain the following characterization of (7): If $x \in S_h$, then

$$\sup_n h(X_n) = \infty, \qquad \mathbf{P}_x^h\text{-a.s.} \quad \iff \quad \inf_n h(X_n) = 0, \qquad \mathbf{P}_x\text{-a.s.}$$

This can also be expressed in terms of the crossing time $\tau_c$ and its dual companion,

$$\tau_c^* := \inf\{n \geq 0 : h(\check{X}_n^*) < U\}$$

as follows,

$$\tau_c < \infty, \qquad \mathbb{P}\text{-a.s.} \quad \iff \quad \tau_c^* < \infty, \qquad \mathbb{P}^*\text{-a.s.}$$

**4. Examples: random walks.** In this section we apply our theorems to random walks $(X_n)$ in $S = \mathbb{R}^d$, that is, Markov chains, whose transition kernels are shift-invariant,

$$P(x, dy) = P(0, -x + dy) \qquad \forall x, y \in \mathbb{R}^d.$$

We assume that the kernel is nondegenerate.

4.1. *Random walk with drift.* First we consider the situation that the random walk has a drift, which means that the increments have finite expectation unequal to zero,

$$\mu := \mathbf{E}_0 X_1 = \int x\, P(0, dx) \neq 0. \tag{9}$$

Furthermore, we assume that

$$\varphi(u) := \int P(0, dx) e^{\langle u, x \rangle} < \infty \tag{10}$$

for all $u \in \mathbb{R}^d$, where $\langle \cdot, \cdot \rangle$ denotes the ordinary dot product in $\mathbb{R}^d$. Note that $\mu = \operatorname{grad} \varphi(0)$.

Then, as it is easy to see, the functions

$$h_u(x) := e^{\langle u, x \rangle}, \qquad u \in C,$$

are harmonic, where

$$C := \{u : \varphi(u) = 1\}.$$

The $h_u$-transformed kernels are given by

$$P^{h_u}(x, dy) = P(x, dy) e^{\langle u, y - x \rangle}.$$

Note that these kernels are again shift-invariant, which means that $(X_n)$ is a random walk also with respect to the measures $\mathbf{P}_x^{h_u}$. The corresponding drift is

$$\mu_u := \mathbf{E}_0^{h_u} X_1 = \int x\, P^{h_u}(0, dx) = \int x e^{\langle u, x \rangle}\, P(0, dx) = \operatorname{grad} \varphi(u).$$

It is well known that $C$ is a smooth manyfold in $\mathbb{R}^d$ of codimension 1. For us the following elementary property is of interest: For $u \in C$, $u \neq 0$, let $\varphi_u(\lambda) := \varphi(\lambda \cdot u) = \mathbf{E}_0 e^{\lambda \langle u, X_1 \rangle}$, $\lambda \in \mathbb{R}$, the moment generating function of $\langle u, X_1 \rangle$. By assumption this random variable is nondegenerate, therefore, $\varphi_u$ is strictly convex. Furthermore, $\varphi_u(\lambda) = 1$ for exactly two values, namely, 0 and 1, therefore, $\varphi_u'(0) < 0$ and $\varphi_u'(1) > 0$. By means of the chain rule, these derivatives are easily calculated as $\langle u, \operatorname{grad} \varphi(0) \rangle$ and $\langle u, \operatorname{grad} \varphi(u) \rangle$, such that

$$\langle u, \mu \rangle < 0, \qquad \langle u, \mu_u \rangle > 0 \qquad \text{for all } u \in C, u \neq 0.$$

Applying the law of large numbers to the random walk $(X_n)$, we obtain

$$h_u(X_n) \to 0, \qquad \mathbf{P}_x\text{-a.s., and} \quad h_u(X_n) \to \infty, \qquad \mathbf{P}_x^{h_u}\text{-a.s.}$$

Thus, (8) is fulfilled. Coming to path decompositions, we obtain that the corresponding splitting and crossing times are a.s. finite,

$$T < \tau_c < \infty \qquad \text{a.s.} \quad \text{and} \quad T^* < \tau_c^* < \infty \qquad \text{a.s.}$$

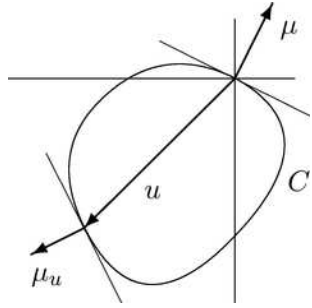

Let us describe our first path decomposition in more detail. We choose $o := 0 \in \mathbb{R}^d$, thus, $h_u(o) = 1$. Let $Y$ be as in theorem 2, that is, $\mathbb{P}(Y > y) = 1/y$ for $y > 1$, and let

$$M := \|u\|^{-1} \log Y.$$

Then $\mathbb{P}(M > m) = \mathbb{P}(Y > e^{\|u\|m}) = e^{-\|u\|m}$, thus, $M$ is exponentially distributed with expectation $\|u\|^{-1}$. Moreover,

$$h_u(x) > Y \quad \Longleftrightarrow \quad \frac{\langle u, x \rangle}{\|u\|} > M \quad \Longleftrightarrow \quad x \in H_u(M),$$

where $H_u(M) := \{x : \langle u, x \rangle > \|u\| M\}$ is a halfspace in $\mathbb{R}^d$, which has distance $M$ from the origin. Altogether we obtain for each $u \in C, u \neq 0$, a different recipy to simulate the original random walk with initial state 0.

1. Let $M$ be exponentially distributed with expectation $\|u\|^{-1}$. Let $\hat{X}$ be an independent random walk with transition kernel $P^{h_u}$.
2. Define

$$T := \sup\{n \geq 0 : X_n \notin H_u(M), \langle u, \hat{X}_m \rangle < \langle u, \hat{X}_n \rangle \, \forall m < n\}.$$

   This is the first moment when $(\hat{X}_n)$ is closest to the halfspace $H_u(M)$, but has not yet entered it.
3. Let $\check{X}$ be a random walk with transition kernel $P$ and initial position $\check{X}_0 = \hat{X}_T$, conditioned not to leave the halfspace, and put

$$\bar{X}_n := \begin{cases} \hat{X}_n, & n < T, \\ \check{X}_{n-T}, & n \geq T. \end{cases}$$

   Then $(\bar{X}_n)$ is a random walk with transition kernel $P$.

In particular, if $(X_n)$ is a one-dimensional random walk with negative drift, then $C$ contains (besides 0) just one positive number $u$. The $h$-transformed process $\hat{X}$ is a random walk with positive drift, and $T$ is the last moment of a record value, before $\hat{X}$ exceeds $M$. The corresponding path decomposition is in complete analogy with Williams' path decomposition for Brownian motion with drift.

The dual decomposition has a similar form. There are plenty of other harmonic functions for multi-dimensional random walks with drift, which lead to further (albeit less transparent) path decompositions.

4.2. *Driftless random walk with absorbtion.* A random walk $(Z_n)$ on the real line $S = \mathbb{R}$ with expected increments

$$\mathbf{E}_0 Z_1 = \int x\, P(0, dx) = 0$$

does not possess any interesting harmonic functions. Therefore we change the setting by making all negative states into absorbing states, that is, we consider a driftless random walk with absorbtion in $(-\infty, 0)$,

$$X_n := Z_{\tau \wedge n} \qquad \text{where } \tau := \inf\{i \geq 0 : Z_i < 0\}.$$

Then $P(x, dy) = P(0, -x + dy)$, if $x \geq 0$, and $P(x, dy) = \delta_x(dy)$, if $x < 0$, where $\delta_x$ is the Dirac-measure in $x$. It is well known [Bertoin and Doney (1994)] that

$$h(x) := \begin{cases} \sum_{k \geq 0} \mathbb{P}(H_k \leq x), & x \geq 0, \\ 0, & x < 0, \end{cases}$$

defines a harmonic function with respect to $P$, where $H_k$ are the strictly descending ladder epochs for the random walk $(Z_n)$ with initial state 0.

The process $(X_n)$ also has under the measure $\mathbf{P}_x^h$, a clear meaning: It is equal in distribution to the process $(Z_n)$, conditioned never to enter $(-\infty, 0)$ [Bertoin and Doney (1994)].

Thus, our theorems offer two dual path decompositions: A path decomposition for *driftless random walks with absorbtion in the negative numbers—at their maximal value*, and a path decomposition for *driftless random walks, conditioned to stay positive—at their minimal value*. Details are left to the reader.

4.3. *Isotropic random walks on homogeneous trees.* Consider a homogeneous tree $T_r$ with $r \geq 3$ edges at each node. The tree can be seen as the Cayley-graph of a free group $G$ generated by $r$ elements $A := \{a_1, \ldots, a_r\}$ with $a_1^2 = \cdots = a_r^2 = e$, where $e$ is the identity. Then each node $x$ corresponds to exactly one formal string

$$x = a_{i_1} a_{i_2} \cdots a_{i_d}$$

with $i_j \neq i_{j+1}$ for all $1 \leq j \leq d-1$. Define the distance of a node $x$ from the center $e$ as

$$|x| := d \qquad \text{if } x \text{ corresponds to the reduced word } x = a_{i_1} a_{i_2} \cdots a_{i_d}$$

and let $c(x,y)$ be the last common vertex in the graph for $x,y \in T_r$. Two reduced nodes $x = a_{i_1}a_{i_2}\cdots a_{i_k}$ and $y = a_{j_1}a_{j_2}\cdots a_{j_l}$ are neighbors if and only if $|k-l| = 1$ and $|y^{-1}x| = 1$, where $y^{-1} = a_{j_l}a_{j_{l-1}}\cdots a_{j_1}$.

Consider an isotropic nearest neighbor random walk on $T_r$, that is, a random walk $(X_n)$ with values in $T_r$ and

$$P(x, xa_i) = \frac{1}{r} \qquad \forall 1 \leq i \leq r.$$

Sawyer and Cartwright have shown [Sawyer (1997) and Cartwright and Sawyer (1991)] that in the isotropic case all minimal harmonic functions of this process are given by

$$h_\omega(x) = (r-1)^{2|c(x,\omega)|-|x|}$$

with $\omega = a_{i_1}a_{i_2}\cdots \in A^\mathbb{N}$ an infinite sequence corresponding to an endpoint at infinity of the tree $T_r$. Thus, the $h_\omega$-transform corresponds to conditioning the isotropic random walk to end in $\omega$.

In this case the constant harmonic function is not minimal, as it is given by

$$1 = \int_{A^\mathbb{N}} h_\omega(x)\mu(d\omega),$$

with $\mu$ being the uniform distribution on $A^\mathbb{N}$ with mass one. This can be easily seen by a bare-hand calculation:

$$\int_{A^\mathbb{N}} h_\omega(x)\mu(d\omega)$$

$$= \frac{1}{r(r-1)^{|x|-1}} \sum_{\omega \in A^{|x|}} (r-1)^{2|c(x,\omega)|-|x|}$$

$$= \frac{1}{r(r-1)^{|x|-1}} \sum_{k=0}^{|x|} \sum_{\substack{\omega \in A^{|x|} \\ |c(x,\omega)|=k}} (r-1)^{2k-|x|}$$

$$= \frac{1}{r(r-1)^{|x|-1}} \Bigg( \underbrace{(r-1)^{|x|}(r-1)^{-|x|}}_{k=0}$$

$$+ \sum_{k=1}^{|x|-1} \underbrace{(r-2)(r-1)^{|x|-k-1}(r-1)^{2k-|x|}}_{k=1,\ldots,|x|-1} + \underbrace{(r-1)^{|x|}}_{k=|x|} \Bigg)$$

$$= \frac{1}{r(r-1)^{|x|-1}}(1 + (r-1)^{|x|-1} - 1 + (r-1)^{|x|})$$

$$= 1.$$

Thus, the isotropic random walk can be constructed by first chosing uniformly an endpoint $\omega$ and then starting a $h_\omega$-transformed random walk.

The transition kernel $P^{h_\omega}$ of such a minimal harmonic function $h_\omega$ then is given by

$$P^{h_\omega}(x,y) = \begin{cases} \dfrac{r-1}{r}, & \text{if } d_\omega(x) < d_\omega(y), \\ \dfrac{1}{r(r-1)}, & \text{if } d_\omega(x) > d_\omega(y), \end{cases}$$

for adjacent nodes $x, y$, with the distance of $x$ to $\omega$ defined as $d_\omega(x) := 2|c(x,\omega)| - |x|$.

The rich harmonic function space allows both the direct decomposition and the dual decomposition. As in the case of the random walk, again we have

$$\inf_n h_\omega(X_n) = 0, \quad \mathbf{P}\text{-a.s.} \quad \text{and} \quad \sup_n h_\omega(X_n) = \infty, \quad \mathbf{P}^{h_\omega}\text{-a.s.},$$

which means for the corresponding splitting times $T$ and $T^*$:

$$T < \tau_c < \infty \quad \text{a.s.} \quad \text{and} \quad T^* < \tau_c^* < \infty \quad \text{a.s.}$$

Fix a harmonic function $h_\omega$. For both decompositions we still need the kernels of the conditioned chains and, therefore, we need

$$q_s(x) := \lambda \mathbf{P}_x(h_\omega(X_i) \leq s \,\forall i \in \mathbb{N}_0),$$
$$q_s^*(x) := \lambda^* \mathbf{P}_x^{h_\omega}(h_\omega(X_i) \geq s \,\forall i \in \mathbb{N}_0)$$

with some $\lambda, \lambda^* > 0$. The kernel $Q_s$ of a conditioned $P$-chain (resp. $Q_s^*$ of a conditioned $P^h$ chain) is then given by

$$Q_s(x,y) = 1/q_s(x) P(x,y) q_s(y), \qquad Q_s^*(x,y) = 1/q_s^*(x) P^h(x,y) q_s^*(y).$$

First we deduce a closed form for $q_s$. The definition of $q_s$ immediately gives us

$$q_s(x) = \mathbf{P}_x((r-1)^{d_\omega(X_i)} \leq s \,\forall i \in \mathbb{N}_0).$$

Note that both kernels $P$ and $P^{h_\omega}$, plus the harmonic functions $h_\omega$, only depend on $d_\omega(x) = 2|c(x,\omega)| - |x|$, and, therefore, the function $q_s$ also only depends on the value of $d_\omega(x)$. Thus, $q_s$ can be reexpressed as

$$q_s(x) = \tilde{q}\left(d_\omega(x) - \left\lfloor \frac{\log s}{\log(r-1)} \right\rfloor\right)$$

with $\tilde{q}: \mathbb{Z} \to \mathbb{R}_+$ defined as

$$\tilde{q}(n) := \lambda \mathbf{P}(d_\omega(X_i) \leq 0 \,\forall i \in \mathbb{N}_0 \mid d_\omega(X_0) = n)$$

with $\lambda > 0$ and the boundary condition $\tilde{q}(0) = 1$. Furthermore, $\tilde{q}$ has to satisfy the following recurrence equation as a direct consequence of the harmonicity of $q_s$ on the set $\{x : h(x) \leq s\}$:

$$\tilde{q}(i) = \frac{1}{r}\tilde{q}(i+1) + \frac{r-1}{r}\tilde{q}(i-1).$$

Solving this recurrence with the two conditions $\tilde{q}(0) = 1$ and $\tilde{q}(i) = 0$ for $i > 0$ leads to the solution

$$\tilde{q}(n) = \begin{cases} \dfrac{(r-1)^n - r + 1}{2 - r}, & \text{if } n \leq 0, \\ 0, & \text{else.} \end{cases}$$

Now we can construct the first decomposition of a free isotropic random walk on $T_r$ as follows (refer to the decomposition of a random walk for more details):

1. Let $M$ be exponentially distributed with expectation 1. Let $\hat{X}$ be an independent random walk on $T_r$ with transition kernel $P^{h_\omega}$.
2. Define

$$T := \sup\{n \geq 0 : d_\omega(X_n) \leq M, d_\omega(X_m) < d_\omega(X_n) \,\forall\, m < n\}.$$

3. Let $\check{X}$ be a random walk with transition kernel $Q_{h_\omega(\hat{X}_T)}$ and initial position $\check{X}_0 = \hat{X}_T$, and put

$$\bar{X}_n := \begin{cases} \hat{X}_n, & n < T, \\ \check{X}_{n-T}, & n \geq T. \end{cases}$$

Then $(\bar{X}_n)$ is an isotropic random walk on $T_r$.

Some similar considerations in the case of a conditioned $P^{h_\omega}$-process lead to the following recurrence equation for $\tilde{q}^*$,

$$\tilde{q}^*(i) = \frac{r-1}{r}\tilde{q}^h(i+1) + \frac{1}{r}\tilde{q}^h(i-1)$$

with the solution given by

$$\tilde{q}^*(i) = \frac{(r-1)^i - 1}{(r-1)^{i-1}}.$$

The dual decomposition of a $P^{h_\omega}$-chain first starts an isotropic random walk on $T_r$ and attaches a conditioned $P^{h_\omega}$-chain with kernel $Q_s^*$. Details are left to the reader.

Fachbereich Mathematik
Universität Frankfurt
Postfach 11 19 32
Frankfurt am Main 60054
Germany
e-mail: kersting@math.uni-frankfurt.de
e-mail: memisogl@math.uni-frankfurt.de